\documentclass{birkjour_t2}
\usepackage{amsmath}
\usepackage{amsfonts}
\usepackage{amssymb}
\usepackage[english]{babel}
\usepackage{color}
\usepackage{graphicx}
\usepackage{lmodern}
\usepackage{amsthm}
\usepackage{mathrsfs}
\usepackage{microtype}
\usepackage{mathscinet}
\usepackage{enumitem}
\usepackage[cal=boondoxo,bb=ams]{mathalfa}
\usepackage{hyperref}
\hypersetup{hidelinks}

\renewcommand{\theequation}{\arabic{section}.\arabic{equation}}

\newtheorem{thm}{Theorem}[section]

\newtheorem{prop}[thm]{Proposition}

\newtheorem{rem}[thm]{Remark}

\DeclareMathOperator\supp{supp}

\DeclareMathOperator\Fuj{Fuj}

\DeclareMathOperator\loc{loc}
\DeclareMathOperator\lin{lin}

\begin{document}

%
%
%
%
%
%
%
%

\title[Weakly coupled system of semilinear scale-invariant wave equations]
 {Weakly coupled system of semilinear wave equations with distinct scale-invariant terms in the linear part}

\author[W. Chen]{Wenhui Chen}
\address{Institute of Applied Analysis, Faculty for Mathematics and Computer Science\\
	 Technical University Bergakademie Freiberg\\
	  Pr\"{u}ferstra{\ss}e 9\\
	   09596 Freiberg\\
	    Germany}
\email{wenhui.chen@student.tu-freiberg.de}

\author[A. Palmieri]{Alessandro Palmieri}

\address{Institute of Applied Analysis, Faculty for Mathematics and Computer Science\\
	 Technical University Bergakademie Freiberg\\
	  Pr\"{u}ferstra{\ss}e 9\\
	  09596 Freiberg\\
	  Germany}
\email{alessandro.palmieri.math@gmail.com}



\date{January 1, 2004}

\begin{abstract}
In this work we determine the critical exponent for a weakly coupled system of semilinear wave equations with distinct scale-invariant lower order terms, when these terms make both equations in some sense ``parabolic-like''. For the blow-up result the test functions method is applied, while for the global existence (in time) results we use $L^2-L^2$ estimates with additional $L^1$ regularity.
\end{abstract}

\maketitle

\section{Introduction}

In this paper we consider the weakly coupled system of semilinear wave equations with scale-invariant damping and mass terms with different multiplicative constants in the lower order terms
\begin{align}\label{weakly coupled system}
\begin{cases}
u_{tt}-\Delta u +\frac{\mu_1}{1+t}u_t +\frac{\nu_1^2}{(1+t)^2}u = |v|^p,  & x\in \mathbb{R}^n, \ t>0,  \\
v_{tt}-\Delta v +\frac{\mu_2}{1+t}v_t +\frac{\nu_2^2}{(1+t)^2}v = |u|^q,  & x\in \mathbb{R}^n, \ t>0, \\
(u,u_t,v,v_t)(0,x)= (u_0,u_1,v_0,v_1)(x), & x\in \mathbb{R}^n,
\end{cases}
\end{align}
where $\mu_1,\mu_2,\nu_1^2,\nu_2^2$ are nonnegative constants and $p,q>1$. 

As in the case of a single semilinear wave equation with scale-invariant damping and mass term, it turns out that the quantities
\begin{align}\label{def deltas}
\delta_j \doteq (\mu_j-1)^2 -4\nu_j^2, \qquad j=1,2,
\end{align} are 
useful to describe some of the properties of the model \eqref{weakly coupled system} as, for example, the critical exponent. 

First we describe the meaning of the critical exponent for a semilinear weakly coupled system. Let us introduce the notations 
\begin{align}\label{def alphas}
\alpha_j \doteq \frac{1}{2}\big(\mu_j+1-\sqrt{\delta_j}\big) \qquad j=1,2.
\end{align}

In the case  $\delta_1,\delta_2\geqslant (n+1)^2$, for \eqref{weakly coupled system} the critical exponent is given by
\begin{equation} \label{critical exponent intro}
E=E(p,q,\alpha_1,\alpha_2)\doteq \max \left\{\frac{p+1}{pq-1}-\frac{\alpha_1-1}{2}, \frac{q+1}{pq-1}-\frac{\alpha_2-1}{2}\right\}=\frac{n}{2},
\end{equation}  that is, if $E<\frac{n}{2}$ (supercritical case), then, there exists a unique global solution for small data; else, if $E\geqslant \frac{n}{2}$ (subcritical or critical case), the local in time solution blows up in finite time.

Although we will be able to determine a blow-up result in the case in which $\delta_1, \delta_1\geqslant 0$, due to the fact that a single scale-invariant wave equation shows properties analogous to those of the classical damped wave equation only for large values of the parameter $\delta$, we will find a sharp result only in the case in which $\delta_1,\delta_2\geqslant (n+1)^2$ (see also \cite{NPR16} for further explanations about this condition).

We recall now some historical background to \eqref{weakly coupled system}. Over the last years, semilinear weakly coupled systems have been widely studied. 

Let us begin with the semilinear weakly coupled system of classical wave equations 
\begin{equation}\label{weaklycoupledwave}
\left\{\begin{aligned}
&u_{tt}-\Delta u=|v|^{p},&\quad &x\in\mathbb{R}^n,\,\,t>0,\\
&v_{tt}-\Delta v=|u|^{q},&\quad &x\in\mathbb{R}^n,\,\,t>0,\\
&(u,u_t,v,v_t)(0,x)=(u_0,u_1,v_0,v_1)(x),&\quad &x\in\mathbb{R}^n,
\end{aligned}\right.
\end{equation}
with $p,q>1$. On the one hand, for the single semilinear wave equation we refer to the classical works \cite{Joh79,Kato80,Str81,Gla81-g,Gla81-b,Sid84,Sch85,LS96,GLS97,Jiao03,YZ06,Zhou07,LaiZhou14}, where the so-called Strauss exponent $p_0(n)$ is proved to be the critical exponent,  $p_0(n)$ being the positive root of the quadratic equation $(n-1)p^2-(n+1)p-2=0$. On the other hand, collecting the results from \cite{DelS97,DGM,DM,AKT00,KT03,Kur05,GTZ06,KTW12}, we find that the critical exponent for \eqref{weaklycoupledwave} is described by the condition
\begin{equation*}
\max\left\{\frac{p+2+q^{-1}}{pq-1},\frac{q+2+p^{-1}}{pq-1}\right\}=\frac{n-1}{2}.
\end{equation*}  

Let us recall some results for the semilinear weakly coupled system of classical damped wave equations 
\begin{equation}\label{weaklycoupleddamped}
\left\{\begin{aligned}
&u_{tt}-\Delta u+u_t=|v|^{p},&\quad &x\in\mathbb{R}^n,\,\,t>0,\\
&v_{tt}-\Delta v+v_t=|u|^{q},&\quad &x\in\mathbb{R}^n,\,\,t>0,\\
&(u,u_t,v,v_t)(0,x)=(u_0,u_1,v_0,v_1)(x),&\quad &x\in\mathbb{R}^n,
\end{aligned}\right.
\end{equation}
with $p,q>1$. For the single semilinear damped wave equation  $p_{\Fuj}(n)\doteq 1+\frac{2}{n}$ is the critical exponent, we refer to the classical works \cite{TodYor01,Zhang01,IkeTa05} for further details.  The critical exponent for \eqref{weaklycoupleddamped} is described by the condition
\begin{equation*}
\max\left\{\frac{p+1}{pq-1},\frac{q+1}{pq-1}\right\}=\frac{n}{2}.
\end{equation*} 
The authors of \cite{SunWang2007} investigated the critical exponent for $n=1,3$. 
In \cite{Narazaki2009} the author generalized the global existence result to $n=1,2,3$ and improved the time decay estimates when $n=3$. After that, in \cite{Nishi12} the asymptotic profile for global solutions has been derived in low dimensions $n=1,2,3$. Then, in \cite{NishiharaWakasugi} global existence and blow-up in finite time results for any space dimension $n$ were determined, where the proof of the global (in time) existence of energy solutions is based on a weighted energy method. Consequently, in \cite{NishiWak15} the previous result has been extended for a semilinear weakly coupled system of $k\geqslant 2$ damped wave equations. In comparison to the critical exponent for \eqref{weaklycoupleddamped}, we observe a translation in the critical exponent for the model that we consider in this work, which is due to the presence of the lower order scale-invariant terms.
We also mention that several generalizations of \eqref{weaklycoupleddamped} are possible in different ways. On the one hand, the weakly coupled system of damped waves with time-dependent coefficients in the dissipation terms is studied, for example, in \cite{NishiWak,Djaouti,Djaouti1}. In particular, in \cite{Djaouti,Djaouti1} the global existence of solutions is proved, when initial data are supposed to belong to different classes of regularity. On the other hand, in \cite{Dabb15} semilinear weakly coupled systems are studied replacing the classical damping terms with structural damping terms.  Finally, in \cite{ChenRei18} a semilinear weakly coupled system of damped elastic waves is studied. In this latter case, the system is coupled not only in the nonlinear terms but also in the linear ones.

Recently, the Cauchy problem 
\begin{align}\label{scale inv eq}
\begin{cases}
u_{tt}-\Delta u +\frac{\mu}{1+t}u_t +\frac{\nu^2}{(1+t)^2}u = |u|^p,  & x\in \mathbb{R}^n, \ t>0,  \\
(u,u_t)(0,x)= (u_0,u_1)(x),& x\in \mathbb{R}^n, 
\end{cases}
\end{align} has attracted a lot of attention, where $\mu,\nu^2$ are nonnegative constants and $p>1$ and, analogously to what we did for the system we denote $\delta\doteq (\mu -1) ^2-4\nu^2$. The value of $\delta$ has a strong influence on some properties of solutions to \eqref{scale inv eq} and to the corresponding homogeneous linear equation. According to \cite{Abb15,Wakasugi14,DLR15,DabbLuc15,Wakasa16,NPR16,PalRei17,Pal17,LTW17,IS17,PR17vs,TL1709,TL1711,Pal18odd,Pal18even,DabbPal18,PT18,KatoSak18,Lai18} for $\delta\geqslant 0$ the model in \eqref{scale inv eq} is somehow an intermediate model between the semilinear free wave equation and the semilinear classical damped equation, whose critical exponent is $p_{\text{Fuj}}(n+\alpha-1)$ for $\delta\geq (n+1)^2$,  where $\alpha$ is defined analogously as in \eqref{def alphas}, and seems reasonably to be $p_0(n+\mu)$ for small values of delta. In this paper we will deal with the system \eqref{weakly coupled system} and we will investigate how the interaction between the powers $p,q$ in the nonlinearities provides either the global in time existence of the solution or the blow-up in finite time.

{\bf Notations:} Throughout this paper we will use the following notations: $B_R$ denotes the ball around the origin with radius $R$; $f\lesssim g$ means that there exists a positive constant $C$ such that $f\leqslant C g$ and, similarly, for $f\gtrsim g$; 
finally, as in the introduction, $p_{\Fuj}(n)$ and $p_0(n)$ denote the Fujita exponent and the Strauss exponent, respectively.
\section{Main results} \label{Section main results}

In \cite{NPR16} a blow-up result is proved for \eqref{scale inv eq} provided that $\delta\geqslant 0$ by using the so-called test function method in the case in which the exponent of the power nonlinearity is smaller than or equal to  $p_{\Fuj}(n+\alpha -1)$. 
In the next result we will generalize that result for the weakly coupled system \eqref{weakly coupled system}. Let us underline that, due to the presence of generally different coefficients in the linear terms of lower order, a new phenomenal appears, that cannot be observed for single equations or for weakly coupled systems with the same linear part (for example, in the case of  \eqref{weakly coupled system} when $\mu_1=\mu_2$ and $\nu_1^2=\nu_2^2$). More precisely, a restriction from below either for $p$ or for $q$ is necessary to get the desired result (see also Remark \ref{Rem lower bound p or lower bound q}).

\begin{thm}[Blow-up result] \label{Thm blow-up MTF} Let $\mu_1,\mu_2,\nu_1^2,\nu_2^2$ be nonnegative constants such that $\delta_1,\delta_2\geqslant 0$ and let $(u_0,u_1,v_0,v_1)\in \big(H^1(\mathbb{R}^n)\times L^2(\mathbb{R}^n)\big)^2$ be initial data such that
	\begin{align}
	\liminf_{R\to\infty} \int_{\mathbb{R}^n}\bigg(\frac{1}{2}\big(\mu_1-1+\sqrt{\delta_1}\,\big)\,u_0(x)+u_1(x)\bigg) dx >0, \label{assumptions u0,u1 blowup} \\
	\liminf_{R\to\infty} \int_{\mathbb{R}^n}\bigg(\frac{1}{2}\big(\mu_2-1+\sqrt{\delta_2}\,\big)\,v_0(x)+v_1(x)\bigg) dx >0. \label{assumptions v0,v1 blowup}
	\end{align}
	If $p,q>1$ satisfy the relations
	\begin{align} \label{condition p,q blowup}
	&\max\bigg\{\frac{p+1}{pq-1}-\frac{\alpha_1}{2}, \frac{q+1}{pq-1}-\frac{\alpha_2}{2}\bigg\} -\frac{n-1}{2}\geqslant 0, \\
	& \mbox{either} \quad p>\frac{1+\alpha_1}{1+\alpha_2} \quad \mbox{or} \quad q>\frac{1+\alpha_2}{1+\alpha_1}, \label{lower bound p or q blowup} 
	\end{align} then, \eqref{weakly coupled system} has no globally in time weak solutions; that is, if $(u,v)\in L^q_{\loc}([0,T)\times \mathbb{R}^n) \times L^p_{\loc}([0,T)\times \mathbb{R}^n)$ is a local in time weak solution with maximal lifespan $T$, then, $T<\infty$.
	
\end{thm}

\begin{rem} \label{Rem lower bound p or lower bound q} We should point out that at least one of the conditions in \eqref{lower bound p or q blowup} is trivially true. Indeed, since $\alpha_j\geqslant 0$ for $j=1,2$, then, in the case $\alpha_1>\alpha_2$, it holds $q>1> \tfrac{1+\alpha_2}{1+\alpha_1}$; while in the case $\alpha_1<\alpha_2$, we have $p>1> \tfrac{1+\alpha_1}{1+\alpha_2}$.
	Moreover, if $\alpha_1=\alpha_2$, for example when we have the same coefficients  for the linear terms in \eqref{weakly coupled system}, then, no additional requirement on $p$ or on $q$ is necessary.
\end{rem}

Due to the blow-up result in Theorem \ref{Thm blow-up MTF}, we expect 
\begin{align}\label{critical exponent weakly coupled system}
E(p,q,\alpha_1,\alpha_2)-\tfrac{n}{2}\equiv \max\big\{F(p,q,n,\alpha_1), F(q,p,n,\alpha_2)\big\}=0
\end{align} to be the critical exponent for the semilinear system \eqref{weakly coupled system} in the case in which both linear parts are somehow ``parabolic-like'' (see Remark \ref{rmk parabolic}), where 
\begin{align*}
F(p,q,n,\alpha)\doteq \tfrac{p+1}{pq-1}-\tfrac{n+\alpha-1}{2}.
\end{align*}

Correspondingly to the case of a single semilinear equation with power nonlinearity, we mean that if $p,q>1$ satisfy 
\begin{align}\label{critical exponent weakly coupled system <0}
\max\big\{F(p,q,n,\alpha_1), F(q,p,n,\alpha_2)\big\}<0,
\end{align}
then, there exists a unique global solution for small initial data, whereas a local in time solution with some integral sign assumptions for the data blows up in finite time if the left hand side in \eqref{critical exponent weakly coupled system} is nonnegative.

As we have already shown a result for the necessity part, now we want to investigate the sufficiency part. Before doing that, we shall clarify under which necessary condition the left hand side in \eqref{critical exponent weakly coupled system} can be negative. For this purpose we introduce the following notations:
\begin{align*}
\widetilde{p}\,(n,\alpha_1,\alpha_2) &\doteq \tfrac{n+\alpha_1+1}{n+\alpha_2-1}=1+\tfrac{2+(\alpha_1-\alpha_2)}{n+\alpha_2-1}, \\
\widetilde{q}\,(n,\alpha_1,\alpha_2) &\doteq \tfrac{n+\alpha_2+1}{n+\alpha_1-1}=1+\tfrac{2+(\alpha_2-\alpha_1)}{n+\alpha_1-1}.
\end{align*}
Let us remark that if $p\leqslant \widetilde{p}\,(n,\alpha_1,\alpha_2)$ and $q\leqslant \widetilde{q}\,(n,\alpha_1,\alpha_2)$, then $F(p,q,n,\alpha_1)\geqslant 0$. Indeed,
\begin{align*}
\big(q-\tfrac{2}{n+\alpha_1-1}\big)\leqslant \tfrac{n+\alpha_2-1}{n+\alpha_1-1} & \qquad \Rightarrow \qquad p\big(q-\tfrac{2}{n+\alpha_1-1}\big)\leqslant \tfrac{n+\alpha_1+1}{n+\alpha_1-1} = 1+\tfrac{2}{n+\alpha_1-1} \\
& \qquad \Rightarrow \qquad pq-1 \leqslant  \tfrac{2}{n+\alpha_1-1} (p+1) \\
& \qquad \Rightarrow \qquad \tfrac{n+\alpha_1-1}{2} \leqslant  \tfrac{p+1}{pq-1}.
\end{align*} Analogously, under the same assumptions on $p$ and $q$ it holds $F(q,p,n,\alpha_2)\geqslant 0$. Summarizing, $p\leqslant \widetilde{p}\,(n,\alpha_1,\alpha_2)$ and $q\leqslant \widetilde{q}\,(n,\alpha_1,\alpha_2)$ imply that the left hand side in \eqref{critical exponent weakly coupled system} is nonnegative. Therefore,
\begin{align*}
\max\big\{F(p,q,n,\alpha_1), F(q,p,n,\alpha_2)\big\} <0 \ \ \mbox{implies}\ \  p > \widetilde{p}\,(n,\alpha_1,\alpha_2)  \ \ \mbox{or} \ \ q> \widetilde{q}\,(n,\alpha_1,\alpha_2).
\end{align*}

Consequently, in order to prove the global in time existence for small data solutions provided that $(p,q)$ satisfies \eqref{critical exponent weakly coupled system <0}, we may consider separately the following three subcases:
\begin{align}
& p > \widetilde{p}\,(n,\alpha_1,\alpha_2)  \ \ \mbox{and} \ \ q> \widetilde{q}\,(n,\alpha_1,\alpha_2), \label{supercritical p and q} \\
& p \leqslant \widetilde{p}\,(n,\alpha_1,\alpha_2)  \ \ \mbox{and} \ \ q> \widetilde{q}\,(n,\alpha_1,\alpha_2), \label{supercritical q and subcritical p}\\
& p > \widetilde{p}\,(n,\alpha_1,\alpha_2)  \ \ \mbox{and} \ \ q\leqslant \widetilde{q}\,(n,\alpha_1,\alpha_2). \label{supercritical p and subcritical q}
\end{align} 
More precisely, in the case \eqref{supercritical p and q} no loss of decay with respect to the corresponding linear problem will appear in the decay estimates. On the other hand, in the case where $p,q>1$ fulfill \eqref{supercritical q and subcritical p} (respectively \eqref{supercritical p and subcritical q}), because different power source nonlinearities have different influence on conditions for the global (in time) existence of solutions, we allow the effect of the loss of decay.

Before stating these global existence results, we should recall some known results for the family of parameter dependent linear Cauchy problems
\begin{align}\label{scale inv eq lin}
\begin{cases}
u_{tt}-\Delta u +\frac{\mu}{1+t}u_t +\frac{\nu^2}{(1+t)^2}u =0,  & x\in \mathbb{R}^n, \ t>s,  \\
(u,u_t)(s,x)= (u_0,u_1)(x),& x\in \mathbb{R}^n, 
\end{cases}
\end{align} where the initial time $s\geqslant 0$ is considered because of the lack of invariance by time-translation for this linear model with time-dependent coefficients. For the proofs of the next results we refer to \cite[Theorems 4.6 and 4.7]{PalRei17}.

\begin{prop}\label{Prop Lin Estim} Let $\mu>0$ and $\nu^2$ be nonnegative constants such that $\delta>(n+1)^2$. Let us consider $(u_0,u_1)\in \big(H^1(\mathbb{R}^n)\cap L^1(\mathbb{R}^n)\big)\times \big(L^2(\mathbb{R}^n)\cap L^1(\mathbb{R}^n)\big)$. Then, for all $\kappa\in[0,1]$ the energy solution $u=u(t,x)$ to \eqref{scale inv eq lin} with $s=0$ satisfies the decay estimate
	\begin{align}\label{Lin Estim} \|u(t,\cdot)\|_{\dot{H}^\kappa(\mathbb{R}^n)}&\lesssim \left(\| u_0\|_{H^1(\mathbb{R}^n)\cap L^1(\mathbb{R}^n)}+\| u_1\|_{L^2(\mathbb{R}^n)\cap L^1(\mathbb{R}^n)}\right) (1+t)^{-\kappa-\frac{n}{2}-\alpha+1}.  
	\end{align}
	Moreover, $\| u_t(t,\cdot)\|_{L^2(\mathbb{R}^n)}$ satisfies the same decay estimates as $\|\nabla u(t,\cdot)\|_{L^2(\mathbb{R}^n)}$. 
\end{prop}

\begin{prop}\label{Prop Lin Estim v_0=0 for t>s}
	Let $\mu>0$ and $\nu^2$ be nonnegative constants such that $\delta>(n+1)^2$. Let us assume $u_0=0$ and $ u_1\in L^2(\mathbb{R}^n)\cap L^1(\mathbb{R}^n)$. Then, the energy solution $u=u(t,x)$ to \eqref{scale inv eq lin} satisfies for $t \geq s$ and $\kappa\in[0,1]$ the following estimate
	\begin{equation}\begin{split}\label{Lin Estim v_0=0 for t>s}
	\|u(t,\cdot)\|_{\dot{H}^\kappa(\mathbb{R}^n)}&\lesssim \left(\|u_1\|_{L^1(\mathbb{R}^n)}+(1+s)^{\frac{n}{2}}\|u_1\|_{L^2(\mathbb{R}^n)}\right)(1+s)^{\alpha} (1+t)^{-\kappa-\frac{n}{2}-\alpha+1}. 
	\end{split}
	\end{equation}
	Moreover, $\| u_t(t,\cdot)\|_{L^2(\mathbb{R}^n)}$ satisfies the same decay estimates as $\|\nabla u(t,\cdot)\|_{L^2(\mathbb{R}^n)}$. 
\end{prop}

\begin{rem} \label{rmk parabolic} In the case $\delta>(n+1)^2$ for the linear Cauchy problem \eqref{scale inv eq lin} we get for the $L^2$ norms of the derivatives a better decay rate than the one for the $L^2$ norm of the solution. Considering larger values of $\delta$ we can observe this phenomenon even for derivatives of higher order, with an improved decay rate as well. In this sense, we say that \eqref{scale inv eq lin} is ``parabolic-like'' for large values of $\delta$.
\end{rem}

Now we can state the main global existence results for \eqref{weakly coupled system}. As in the previous propositions, we will work with initial data in the classical energy space with additional $L^1$ regularity, so that the space for the Cauchy data is
\begin{align*}
\mathcal{A}\doteq (H^1(\mathbb{R}^n)\cap L^1(\mathbb{R}^n)\big)\times \big(L^2(\mathbb{R}^n)\cap L^1(\mathbb{R}^n)\big).
\end{align*} 

Let us begin with the subcase \eqref{supercritical p and q}.

\begin{thm}\label{Thm global existence supercritical p,q} Let $\mu_1,\mu_2>1$, $\nu_1^2,\nu_2^2$ be nonnegative constants such that $\delta_1,\delta_2>(n+1)^2$. Let us assume $p,q>1$, satisfying $2\leqslant p,q $ and $p,q \leqslant  \tfrac{n}{n-2}$ if $n\geqslant 3$, such that
	\begin{align*}
	& p > \widetilde{p}\,(n,\alpha_1,\alpha_2)  \ \ \mbox{and} \ \ q> \widetilde{q}\,(n,\alpha_1,\alpha_2).
	\end{align*} Then, there exists a constant $\varepsilon_0>0$ such that for any $(u_0,u_1)\in\mathcal{A}$ and $(v_0,v_1)\in \mathcal{A}$ with $$\|(u_0,u_1)\|_{\mathcal{A}}+\|(v_0,v_1)\|_{\mathcal{A}}\leqslant \varepsilon_0$$ there is a uniquely determined energy solution $(u,v)\in \big(\mathcal{C}([0,\infty),H^1(\mathbb{R}^n))\cap\mathcal{C}^1([0,\infty),L^2(\mathbb{R}^n))\big)^2$ to \eqref{weakly coupled system}.
	Furthermore, the solution $(u,v)$ satisfies the following decay estimates: 
	\begin{align*}
	\|(\nabla u,u_t)(t,\cdot)\|_{L^2(\mathbb{R}^n)} & \lesssim (1+t)^{-\frac{n}{2}-\alpha_1} \big(\|(u_0,u_1)\|_{\mathcal{A}}+\|(v_0,v_1)\|_{\mathcal{A}}\big), \\
	\|u(t,\cdot)\|_{L^2(\mathbb{R}^n)} & \lesssim (1+t)^{-\frac{n}{2}-\alpha_1+1} \big(\|(u_0,u_1)\|_{\mathcal{A}}+\|(v_0,v_1)\|_{\mathcal{A}}\big), \\
	\|(\nabla v,v_t)(t,\cdot)\|_{L^2(\mathbb{R}^n)} & \lesssim (1+t)^{-\frac{n}{2}-\alpha_2} \big(\|(u_0,u_1)\|_{\mathcal{A}}+\|(v_0,v_1)\|_{\mathcal{A}}\big), \\
	\|v(t,\cdot)\|_{L^2(\mathbb{R}^n)} & \lesssim (1+t)^{-\frac{n}{2}-\alpha_2+1} \big(\|(u_0,u_1)\|_{\mathcal{A}}+\|(v_0,v_1)\|_{\mathcal{A}}\big). 
	\end{align*}
\end{thm} 

\begin{rem} If $p,q>1$ satisfy \eqref{supercritical p and q}, then, $F(p,q,n,\alpha_1)<0$ and $F(q,p,n,\alpha_2)<0$. Also, in particular, \eqref{critical exponent weakly coupled system <0} holds. Indeed,
	\begin{align*}
	\big(q-\tfrac{2}{n+\alpha_1-1}\big)>\tfrac{n+\alpha_2-1}{n+\alpha_1-1} & \qquad \Rightarrow \qquad p\big(q-\tfrac{2}{n+\alpha_1-1}\big)>  1+\tfrac{2}{n+\alpha_1-1} \\
	& \qquad \Rightarrow \qquad F(p,q,n,\alpha_1)<0,
	\end{align*} and, similarly, $F(q,p,n,\alpha_2)<0$.
\end{rem}

Let us consider now the subcase \eqref{supercritical q and subcritical p}.

\begin{thm}\label{Thm global existence subcritical p and supercritical q} Let $\mu_1,\mu_2>1$, $\nu_1^2,\nu_2^2$ be nonnegative constants such that $\delta_1,\delta_2>(n+1)^2$. Let us assume $p,q>1$, satisfying $2\leqslant p,q $ and $p,q \leqslant  \tfrac{n}{n-2}$ if $n\geqslant 3$, such that
	\begin{align}
	& p \leqslant \widetilde{p}\,(n,\alpha_1,\alpha_2)  \ \ \mbox{and} \ \ q> \widetilde{q}\,(n,\alpha_1,\alpha_2), \notag \\
	& F(q,p,n,\alpha_2)\equiv\tfrac{q+1}{pq-1}-\tfrac{n+\alpha_2-1}{2}<0. \label{F(q,p,alpha2)<0}
	\end{align} Then, there exists a constant $\varepsilon_0>0$ such that for any $(u_0,u_1)\in\mathcal{A}$ and $(v_0,v_1)\in \mathcal{A}$ with $$\|(u_0,u_1)\|_{\mathcal{A}}+\|(v_0,v_1)\|_{\mathcal{A}}\leqslant \varepsilon_0$$ there is a uniquely determined energy solution $(u,v)\in \big(\mathcal{C}([0,\infty),H^1(\mathbb{R}^n))\cap\mathcal{C}^1([0,\infty),L^2(\mathbb{R}^n))\big)^2$ to \eqref{weakly coupled system}.
	Furthermore, the solution $(u,v)$ satisfies the following estimates:
	\begin{align*}
	\|(\nabla u,u_t)(t,\cdot)\|_{L^2(\mathbb{R}^n)} & \lesssim (1+t)^{-\frac{n}{2}-\alpha_1+\gamma} \big(\|(u_0,u_1)\|_{\mathcal{A}}+\|(v_0,v_1)\|_{\mathcal{A}}\big), \\
	\|u(t,\cdot)\|_{L^2(\mathbb{R}^n)} & \lesssim (1+t)^{-\frac{n}{2}-\alpha_1+1+\gamma} \big(\|(u_0,u_1)\|_{\mathcal{A}}+\|(v_0,v_1)\|_{\mathcal{A}}\big), \\
	\|(\nabla v,v_t)(t,\cdot)\|_{L^2(\mathbb{R}^n)} & \lesssim (1+t)^{-\frac{n}{2}-\alpha_2} \big(\|(u_0,u_1)\|_{\mathcal{A}}+\|(v_0,v_1)\|_{\mathcal{A}}\big), \\
	\|v(t,\cdot)\|_{L^2(\mathbb{R}^n)} & \lesssim (1+t)^{-\frac{n}{2}-\alpha_2+1} \big(\|(u_0,u_1)\|_{\mathcal{A}}+\|(v_0,v_1)\|_{\mathcal{A}}\big),
	\end{align*} where $$0<\gamma=\gamma(p,n,\alpha_1,\alpha_2)\doteq \begin{cases}(n+\alpha_2-1)(\widetilde{p}\,(n,\alpha_1,\alpha_2)-p) & \mbox{if} \ \ p < \widetilde{p}\,(n,\alpha_1,\alpha_2), \\ \epsilon & \mbox{if} \ \ p = \widetilde{p}\,(n,\alpha_1,\alpha_2), \end{cases}$$ represents the loss of decay in comparison with the corresponding decay estimates for the solution $u$ to the linear Cauchy problem with vanishing right hand side (cf. Proposition \ref{Prop Lin Estim}), $\epsilon>0$ being an arbitrarily small constant in the limit case $p=\widetilde{p}\,(n,\alpha_1,\alpha_2)$.
\end{thm}

\begin{rem} Under the assumptions of Theorem \ref{Thm global existence subcritical p and supercritical q}, the condition $F(p,q,n,\alpha_1)<0$ is only apparently not necessary in order to apply a standard contraction argument. Nevertheless, in the moment in which we split the global (in time) existence results in the subcases \eqref{supercritical p and q}, \eqref{supercritical q and subcritical p} and \eqref{supercritical p and subcritical q}, we have already used the condition \eqref{critical exponent weakly coupled system <0} and, thus, the condition $F(p,q,n,\alpha_1)<0$.

\end{rem}

Finally, switching the role of $p$ and $q$ in Theorem \ref{Thm global existence subcritical p and supercritical q}, we get the next result.

\begin{thm}\label{Thm global existence supercritical p and subcritical q}  Let $\mu_1,\mu_2>1$, $\nu_1^2,\nu_2^2$ be nonnegative constants such that $\delta_1,\delta_2>(n+1)^2$. Let us assume $p,q>1$, satisfying $2\leqslant p,q $ and $p,q \leqslant  \tfrac{n}{n-2}$ if $n\geqslant 3$, such that
	\begin{align}
	& p > \widetilde{p}\,(n,\alpha_1,\alpha_2)  \ \ \mbox{and} \ \ q\leqslant \widetilde{q}\,(n,\alpha_1,\alpha_2), \notag \\
	& F(p,q,n,\alpha_1)\equiv\tfrac{p+1}{pq-1}-\tfrac{n+\alpha_1-1}{2}<0. \label{F(p,q,alpha1)<0}
	\end{align} Then, there exists a constant $\varepsilon_0>0$ such that for any $(u_0,u_1)\in\mathcal{A}$ and $(v_0,v_1)\in \mathcal{A}$ with $$\|(u_0,u_1)\|_{\mathcal{A}}+\|(v_0,v_1)\|_{\mathcal{A}}\leqslant \varepsilon_0$$ there is a uniquely determined energy solution $(u,v)\in \big(\mathcal{C}([0,\infty),H^1(\mathbb{R}^n))\cap\mathcal{C}^1([0,\infty),L^2(\mathbb{R}^n))\big)^2$ to \eqref{weakly coupled system}.
	Furthermore, the solution $(u,v)$ satisfies the following estimates:
	\begin{align*}
	\|(\nabla u,u_t)(t,\cdot)\|_{L^2(\mathbb{R}^n)} & \lesssim (1+t)^{-\frac{n}{2}-\alpha_1} \big(\|(u_0,u_1)\|_{\mathcal{A}}+\|(v_0,v_1)\|_{\mathcal{A}}\big), \\
	\|u(t,\cdot)\|_{L^2(\mathbb{R}^n)} & \lesssim (1+t)^{-\frac{n}{2}-\alpha_1+1} \big(\|(u_0,u_1)\|_{\mathcal{A}}+\|(v_0,v_1)\|_{\mathcal{A}}\big), \\
	\|(\nabla v,v_t)(t,\cdot)\|_{L^2(\mathbb{R}^n)} & \lesssim (1+t)^{-\frac{n}{2}-\alpha_2+\bar{\gamma}} \big(\|(u_0,u_1)\|_{\mathcal{A}}+\|(v_0,v_1)\|_{\mathcal{A}}\big), \\
	\|v(t,\cdot)\|_{L^2(\mathbb{R}^n)} & \lesssim (1+t)^{-\frac{n}{2}-\alpha_2+1+\bar{\gamma}} \big(\|(u_0,u_1)\|_{\mathcal{A}}+\|(v_0,v_1)\|_{\mathcal{A}}\big),
	\end{align*} where $$0<\bar{\gamma}=\bar{\gamma}(q,n,\alpha_1,\alpha_2)\doteq \begin{cases}(n+\alpha_1-1)(\widetilde{q}\,(n,\alpha_1,\alpha_2)-q) & \mbox{if} \ \ q < \widetilde{q}\,(n,\alpha_1,\alpha_2), \\ \epsilon & \mbox{if} \ \ q = \widetilde{q}\,(n,\alpha_1,\alpha_2), \end{cases}$$ represents the loss of decay in comparison with the corresponding decay estimates for the solution $v$ to the linear Cauchy problem with vanishing right hand side (cf. Proposition \ref{Prop Lin Estim}), $\epsilon>0$ being an arbitrarily small constant in the limit case $q=\widetilde{q}\,(n,\alpha_1,\alpha_2)$.
\end{thm}

\begin{rem} Also in this case, the condition $F(q,p,n,\alpha_2)<0$ is implicitly used in the previous theorem.
\end{rem}


\section{Blow-up result: Proof of Theorem \ref{Thm blow-up MTF}}

In this section we will employ the so-called \emph{test functions method} (see for example\cite{MP98,MP99,EGKP00,MP01book,MP01h,MP01ep,Zhang01,GMP15}).

Let us assume by contradiction that $(u,v)\in L^q_{\loc}([0,T)\times \mathbb{R}^n) \times L^p_{\loc}([0,T)\times \mathbb{R}^n)$ is a global (in time) weak solution to \eqref{weakly coupled system}, that is $T=\infty$. This means that the integral equalities 
%
%
\begin{align} 
& \iint_{[0,T)\times \mathbb{R}^n} \bigg(\partial_t^2 \psi_1(t,x)-\Delta \psi_1(t,x)-\partial_t \Big(\tfrac{\mu_1}{1+t}\psi_1(t,x)\Big)+\tfrac{\nu_1^2}{(1+t)^2}\psi_1(t,x)\bigg) u(t,x) \,d(t,x)  \notag \\ 
\quad &= \int_{\mathbb{R}^n} \big(\psi_1(0,x)(u_1(x)+\mu_1u_0(x))-\partial_t \psi_1(0,x)u_0(x)\big) dx+\iint_{[0,T)\times \mathbb{R}^n} \psi_1(t,x) |v(t,x)|^p d(t,x),  \label{def weak sol for the system u eq} \\
& \iint_{[0,T)\times \mathbb{R}^n} \bigg(\partial_t^2 \psi_2(t,x)-\Delta \psi_2(t,x)-\partial_t \Big(\tfrac{\mu_2}{1+t}\psi_2(t,x)\Big)+\tfrac{\nu_2^2}{(1+t)^2}\psi_2(t,x)\bigg) v(t,x) \,d(t,x)  \notag \\ 
\quad &= \int_{\mathbb{R}^n} \big(\psi_2(0,x)(v_1(x)+\mu_2 v_0(x))-\partial_t \psi_2(0,x)v_0(x)\big) dx +\iint_{[0,T)\times \mathbb{R}^n} \psi_2(t,x) |u(t,x)|^q d(t,x),  \label{def weak sol for the system v eq}
\end{align}
 are fulfilled for any $(\psi_1,\psi_2)\in \big(\mathcal{C}^\infty_0([0,T)\times \mathbb{R}^n)\big)^2$.

Multiplying the first and the second equation in \eqref{weakly coupled system} by time-dependent functions $g_1=g_1(t)$ and $g_2=g_2(t)$, respectively, we obtain
\begin{align*}
& \partial_t^2(g_1 u) -\Delta (g_1 u)+ \partial_t \bigg(-2g_1' u+ \frac{\mu_1}{1+t}g_1 u\bigg) + \bigg(g_1'' -\frac{\mu_1}{1+t}g_1'+\frac{\mu_1+\nu_1^2}{(1+t)^2}g_1 \bigg) u = g_1 |v|^p, \\
& \partial_t^2(g_2 v) -\Delta (g_2 v)+ \partial_t \bigg(-2g_2' v+ \frac{\mu_2}{1+t}g_2 v\bigg) + \bigg(g_2'' -\frac{\mu_2}{1+t}g_2'+\frac{\mu_2+\nu_2^2}{(1+t)^2}g_2 \bigg) v = g_2 |u|^q.
\end{align*}
If we choose 
\begin{align}\label{def g functions}
g_j(t)\doteq (1+t)^{\alpha_j}, \qquad j=1,2,
\end{align} then, $g_1$ and $g_2$ satisfy
\begin{align*}
g_j'' -\frac{\mu_j}{1+t}g_j'+\frac{\mu_j+\nu_j^2}{(1+t)^2}g_j=0, \qquad j=1,2.
\end{align*}
Therefore, the previous two relations can be written in the divergence form as follows:
\begin{align*}
& \partial_t^2(g_1 u) -\Delta (g_1 u)+ \partial_t \bigg(-2g_1' u+ \frac{\mu_1}{1+t}g_1 u\bigg) = g_1 |v|^p,  \\
& \partial_t^2(g_2 v) -\Delta (g_2 v)+ \partial_t \bigg(-2g_2' v+ \frac{\mu_2}{1+t}g_2 v\bigg)  = g_2 |u|^q. 
\end{align*}

Let us introduce now two bump functions $\eta \in \mathcal{C}^\infty_0([0,\infty)), \phi\in \mathcal{C}^\infty_0(\mathbb{R}^n)$ such that
\begin{itemize}
	\item $\eta$ is decreasing, $\eta=1$ on $[0,\frac{1}{2}]$ and $\supp\eta\subset [0,1]$;
	\item $\phi$ is radial symmetric and decreasing with respect to $|x|$, $\phi=1$ on $B_\frac{1}{2}$ and $\supp\phi\subset B_1$.
\end{itemize}
These functions satisfy the estimates
\begin{align*}
|\eta'(t)| \lesssim \eta(t)^{\frac{1}{r}}, \ |\eta''(t)| \lesssim \eta(t)^{\frac{1}{r}}, \ |\Delta \phi(x)| \lesssim \phi(x)^{\frac{1}{r}}
\end{align*} for any $r>1$ (see \cite{NPR16}, for example). Moreover, since $0\leqslant \eta(t), \phi(x) \leqslant 1$, then, $\eta(t)\leqslant \eta(t)^{\frac{1}{r}}$ and $\phi(x)\leqslant \phi(x)^{\frac{1}{r}}$ for any $r>1$. In particular, we will use these conditions for $r=p,q$. 

Given two positive parameters $\tau$ and $R$, we define
\begin{align*}
\psi_{\tau,R}(t,x)\doteq \eta_\tau(t) \phi_R(x) \quad \mbox{with} \ \  \eta_\tau(t)\doteq \eta\big(\tfrac{t}{\tau}\big) \ \mbox{and} \ \phi_R(x)\doteq \phi\big(\tfrac{x}{R}\big).
\end{align*} Furthermore, we introduce the following two integrals depending on the parameters $\tau, R$:
\begin{align*}
I_{\tau,R} & \doteq \iint_{[0,T)\times \mathbb{R}^n} g_1(t) \, \psi_{\tau,R}(t,x) \,  |v(t,x)|^p d(t,x) = \int_0^\tau \int_{B_R} g_1(t) \, \psi_{\tau,R}(t,x) \, |v(t,x)|^p dx \, dt  , \\ 
J_{\tau,R} & \doteq \iint_{[0,T)\times \mathbb{R}^n} g_2(t) \, \psi_{\tau,R}(t,x) \, |u(t,x)|^q d(t,x) = \int_0^\tau \int_{B_R} g_2(t) \, \psi_{\tau,R}(t,x) \,  |u(t,x)|^q dx \, dt .
\end{align*} 

Applying the integral relation \eqref{def weak sol for the system u eq} to $g_1 \psi_{\tau,R}$, we get
\begin{align*}
I_{\tau,R}& = -\int_{B_R(0)}\Big(g_1(0)(u_1(x)+\mu_1u_0(x))-g_1'(0)u_0(x)\Big)\phi_R(x)\, dx  \\&\quad+\int_0^\tau\int_{B_R(0)} g_1(t)\,u(t,x)\,\partial_t^2\psi_{\tau,R}(t,x)\,dxdt \\&\quad+\int_0^\tau\int_{B_R(0)} \left(2g_1'(t)-\tfrac{\mu_1}{1+t}g_1(t)\right) u(t,x)\partial_t\psi_{\tau,R}(t,x)\,dxdt \\&\quad-\int_0^\tau\int_{B_R(0)} g_1(t)u(t,x)\Delta\psi_{\tau,R}(t,x)\,dxdt\\&\quad+\int_0^\tau\int_{B_R(0)} \Big(\underbrace{g_1''(t)-\tfrac{\mu_1}{1+t}g_1'(t) +\tfrac{\mu_1+\mu_2^2}{(1+t)^2}g_1(t)}_{{}=0}\Big) u(t,x)\psi_{\tau,R}(t,x)\, dxdt \\
&=-\int_{B_R(0)}\left(u_1(x)+\left(\tfrac{\mu_1}{2}-\tfrac{1}{2}+\tfrac{\sqrt{\delta_1}}{2}\right)u_0(x)\right)\phi_R(x)\,dx +\int_\frac{\tau}{2}^\tau\int_{B_R(0)} g_1(t)\,u(t,x)\,\partial_t^2\psi_{\tau,R}(t,x)\,dxdt \\
&\quad +\int_\frac{\tau}{2}^\tau\int_{B_R(0)} \left(2g_1'(t)-\tfrac{\mu_1}{1+t}g_1(t)\right) u(t,x)\partial_t\psi_{\tau,R}(t,x)\,dxdt\\
&\quad -\int_0^\tau\int_{B_R(0)\setminus B_{R/2}(0)} g_1(t)\,u(t,x)\,\Delta\psi_{\tau,R}(t,x)\,dxdt\\
& \doteq-\int_{B_R(0)}\left(u_1(x)+\left(\tfrac{\mu_1}{2}-\tfrac{1}{2}+\tfrac{\sqrt{\delta_1}}{2}\right)u_0(x)\right)\phi_R(x)\,dx+K_1+K_2+K_3.
\end{align*} Let us underline that in the previous chain of equalities we used \begin{align*}
\supp(\partial_t\psi_{\tau,R}),\ \supp(\partial_t^2\psi_{\tau,R})\subset \left[\tfrac{\tau}{2},\tau\right]\times B_R(0) \ \ \mbox{and} \ \ \supp(\Delta\psi_{\tau,R})\subset [0,\tau]\times (B_R(0)\setminus B_{R/2}(0)).
\end{align*}

Thanks to \eqref{assumptions u0,u1 blowup} and to the properties of $\phi_R$, there exists $R_0$ such that for any $R\geq R_0$
\begin{align*}
\int_{B_R(0)}\left(u_1(x)+\left(\tfrac{\mu_1}{2}-\tfrac{1}{2}+\tfrac{\sqrt{\delta_1}}{2}\right)u_0(x)\right)\phi_R(x)\,dx >0.
\end{align*} 
Thus, for $R\geq R_0$ it holds
\begin{align*}
I_{\tau,R}< K_1+K_2+K_3.
\end{align*} Let us separately estimate the integrals $K_1,K_2,K_3$. Let us begin with $K_1$. Since
\begin{align*}
K_1 & = \tau^{-2} \int_\frac{\tau}{2}^\tau\int_{B_R(0)} g_1(t)\,u(t,x)\,\eta''\big(\tfrac{t}{\tau}\big)\, \phi\big(\tfrac{x}{R}\big)\,dxdt \\
& = \tau^{-2} \int_\frac{\tau}{2}^\tau\int_{B_R(0)} g_2(t)^\frac{1}{q}\,u(t,x)\,\eta''\big(\tfrac{t}{\tau}\big)\, \phi\big(\tfrac{x}{R}\big)\, g_1(t) g_2(t)^{\frac{1}{q'}-1}\,dxdt,
\end{align*} where $q'$ denotes the H\"older conjugate of $q$, by H\"older's inequality it follows
\begin{align*}
|K_1| & \leqslant \tau^{-2}\bigg( \int_\frac{\tau}{2}^\tau\int_{B_R(0)} g_2(t)\, |u(t,x)|^q \,|\eta''\big(\tfrac{t}{\tau}\big)|^q\, \phi\big(\tfrac{x}{R}\big)^q\,dxdt \bigg)^{\frac{1}{q}} \bigg( \int_\frac{\tau}{2}^\tau\int_{B_R(0)}  g_1(t)^{q'} g_2(t)^{1-q'}\,dxdt \bigg)^{\frac{1}{q'}} \\
& \lesssim \tau^{-2}\bigg( \int_\frac{\tau}{2}^\tau\int_{B_R(0)} g_2(t)\, |u(t,x)|^q \, \psi_{\tau,R}(t,x)\,dxdt \bigg)^{\frac{1}{q}} \bigg( \int_\frac{\tau}{2}^\tau\int_{B_R(0)} (1+t)^{(\alpha_1-\alpha_2)q'+\alpha_2}\,dxdt \bigg)^{\frac{1}{q'}}.
\end{align*}

If we introduce the parameter dependent integral
\begin{align*}
\widehat{J}_{\tau,R} & \doteq \int_\frac{\tau}{2}^\tau \int_{B_R} g_2(t) \, \psi_{\tau,R}(t,x)\,  |u(t,x)|^q dx \, dt,
\end{align*} then for $\tau>1$ we get from the last inequality
\begin{align*}
|K_1| & \lesssim \tau^{-2+(\alpha_1-\alpha_2)+\frac{\alpha_2+1}{q'}} R^{\frac{n}{q'}}\widehat{J}_{\tau,R}^{\frac{1}{q}}.
\end{align*}

Let us consider now $K_2$. The relation $$2g'_1(t)-\frac{\mu_1}{1+t}g_1(t)=(2\alpha_1-\mu_1)(1+t)^{\alpha_1-1}=\big(1-\sqrt{\delta_1}\big) \,g_1(t)\, (1+t)^{-1}$$ implies
\begin{align*}
K_2 & =\tau^{-1} \int_\frac{\tau}{2}^\tau\int_{B_R(0)} \left(2g_1'(t)-\tfrac{\mu_1}{1+t}g_1(t)\right) u(t,x)\, \eta'\big(\tfrac{t}{\tau}\big)\, \phi\big(\tfrac{x}{R}\big) \,dxdt \\
& =\big(1-\sqrt{\delta_1}\big)\, \tau^{-1} \int_\frac{\tau}{2}^\tau\int_{B_R(0)} g_2(t)^{\frac{1}{q}} u(t,x)\, \eta'\big(\tfrac{t}{\tau}\big)\, \phi\big(\tfrac{x}{R}\big)\, (1+t)^{-1} g_1(t)\, g_2(t)^{\frac{1}{q'}-1}\,dxdt.
\end{align*} Hence, using H\"older's inequality, we arrive at
\begin{align*}
|K_2| & \lesssim \tau^{-1}\bigg( \int_\frac{\tau}{2}^\tau\int_{B_R(0)} g_2(t)\, |u(t,x)|^q \,|\eta'\big(\tfrac{t}{\tau}\big)|^q\, \phi\big(\tfrac{x}{R}\big)^q\,dxdt \bigg)^{\frac{1}{q}} \\ & \qquad \qquad \qquad \times \bigg( \int_\frac{\tau}{2}^\tau\int_{B_R(0)}  (1+t)^{-q'}g_1(t)^{q'} g_2(t)^{1-q'}\,dxdt \bigg)^{\frac{1}{q'}} \\
& \lesssim \tau^{-1}\bigg( \int_\frac{\tau}{2}^\tau\int_{B_R(0)} g_2(t)\, |u(t,x)|^q \, \psi_{\tau,R}(t,x)\,dxdt \bigg)^{\frac{1}{q}} \bigg( \int_\frac{\tau}{2}^\tau\int_{B_R(0)} (1+t)^{(\alpha_1-\alpha_2-1)q'+\alpha_2}\,dxdt \bigg)^{\frac{1}{q'}} \\
& \lesssim \tau^{-2+(\alpha_1-\alpha_2)+\frac{\alpha_2+1}{q'}} R^{\frac{n}{q'}}\widehat{J}_{\tau,R}^{\frac{1}{q}},
\end{align*} for $\tau >1$. Finally, we estimate $K_3$. Applying again H\"older's inequality to 
\begin{align*}
K_3 & = -R^{-2}\int_0^\tau\int_{B_R(0)\setminus B_{R/2}(0)} g_1(t)\,u(t,x)\,\eta\big(\tfrac{t}{\tau}\big)\,\Delta\phi\big(\tfrac{x}{R}\big)\,dxdt \\
& = -R^{-2}\int_0^\tau\int_{B_R(0)\setminus B_{R/2}(0)} g_2(t)^{\frac{1}{q}}\,u(t,x)\,\eta\big(\tfrac{t}{\tau}\big)\,\Delta\phi\big(\tfrac{x}{R}\big)\,  g_1(t) \, g_2(t)^{\frac{1}{q'}-1}\,dxdt,
\end{align*} we find 
\begin{align*}
|K_3| &\leqslant R^{-2}\bigg(\int_0^\tau\int_{B_R(0)\setminus B_{R/2}(0)} g_2(t)\,|u(t,x)|^q\,\eta\big(\tfrac{t}{\tau}\big)^q\,|\Delta\phi\big(\tfrac{x}{R}\big)|^q \,dxdt\bigg)^{\frac{1}{q}}  \\ & \qquad \qquad \qquad \times \bigg(\int_0^\tau\int_{B_R(0)\setminus B_{R/2}(0)}  g_1(t)^{q'} \, g_2(t)^{1-q'}\,dxdt\bigg)^{\frac{1}{q'}} \\
& \lesssim R^{-2}\bigg(\int_0^\tau\int_{B_R(0)\setminus B_{R/2}(0)} g_2(t)\,|u(t,x)|^q\, \psi_{\tau,R}(t,x) \,dxdt\bigg)^{\frac{1}{q}} \\ & \qquad \qquad \qquad \times  \bigg(\int_0^\tau\int_{B_R(0)\setminus B_{R/2}(0)}  (1+t)^{(\alpha_1-\alpha_2)q'+\alpha_2}\,dxdt\bigg)^{\frac{1}{q'}} \\
& \lesssim R^{-2+\frac{n}{q'}}\widetilde{J}_{\tau,R}^{\frac{1}{q}} \ \bigg(\int_0^\tau (1+t)^{(\alpha_1-\alpha_2)q'+\alpha_2}\,dt\bigg)^{\frac{1}{q'}},
\end{align*} where $\widetilde{J}_{\tau,R}$ is given by
\begin{align*}
\widetilde{J}_{\tau,R}\doteq \int_0^\tau\int_{B_R(0)\setminus B_{R/2}(0)} g_2(t) \, \psi_{\tau,R}(t,x)\,  |u(t,x)|^q dx \, dt.
\end{align*} Differently from the estimates for the terms $K_1,K_2$ in this case we need to consider three subcases for the estimate of the $t$-integral on the right hand side of the last inequality for $|K_3|$, because the integral is no longer over $[\frac{\tau}{2},\tau]$ rather on $[0,\tau]$. Indeed, for $\tau>1$ it holds
\begin{align*}
\int_0^\tau (1+t)^{(\alpha_1-\alpha_2)q'+\alpha_2}\,dt \lesssim 
\begin{cases}
\tau^{(\alpha_1-\alpha_2)q'+\alpha_2+1}  & \mbox{if} \ \  (\alpha_1-\alpha_2)q'+\alpha_2 >-1, \\
\log(1+\tau)  & \mbox{if} \ \  (\alpha_1-\alpha_2)q'+\alpha_2 =-1, \\
1  & \mbox{if} \ \  (\alpha_1-\alpha_2)q'+\alpha_2 <-1.
\end{cases}
\end{align*}
Due to \eqref{lower bound p or q blowup}, we are necessary in the first of the  previous cases, since $q>\tfrac{1+\alpha_2}{1+\alpha_1}$ is equivalent to $(\alpha_1-\alpha_2)q'+\alpha_2 >-1$.
Thus,
\begin{align*}
|K_3| & \lesssim \tau^{(\alpha_1-\alpha_2)+\frac{\alpha_2+1}{q'}}  R^{-2+\frac{n}{q'}} \widetilde{J}_{\tau,R}^{\frac{1}{q}} .
\end{align*}

Consequently, combining the previously obtained estimates for $K_1,K_2,K_3$,  we get
\begin{align}\label{estimate I_R}
I_R \lesssim R^{-2+(\alpha_1-\alpha_2)+\frac{n+\alpha_2+1}{q'}} \Big(\widehat{J}_R^{\frac{1}{q}}+\widetilde{J}_R^{\frac{1}{q}}\Big)
\end{align} for $\tau=R>\max\{R_0,1\}$, where for the sake of simplicity of notation we get rid of the second parameter in the subscript in $I, \widehat{J}, \widetilde{J}$.

Applying the integral relation \eqref{def weak sol for the system v eq} to $g_2 \psi_{\tau,R}$, similarly as for the computations for $I_{\tau,R}$, we get 
\begin{align*}
J_{\tau,R}& = -\int_{B_R(0)}\left(v_1(x)+\left(\tfrac{\mu_2}{2}-\tfrac{1}{2}+\tfrac{\sqrt{\delta_2}}{2}\right)v_0(x)\right)\phi_R(x)\,dx +\int_\frac{\tau}{2}^\tau\int_{B_R(0)} g_2(t)\,v(t,x)\,\partial_t^2\psi_{\tau,R}(t,x)\,dxdt \\
&\quad +\int_\frac{\tau}{2}^\tau\int_{B_R(0)} \left(2g_2'(t)-\tfrac{\mu_2}{1+t}g_2(t)\right) v(t,x)\partial_t\psi_{\tau,R}(t,x)\,dxdt\\
&\quad -\int_0^\tau\int_{B_R(0)\setminus B_{R/2}(0)} g_2(t)\,v(t,x)\,\Delta\psi_{\tau,R}(t,x)\,dxdt\\
& \doteq-\int_{B_R(0)}\left(v_1(x)+\left(\tfrac{\mu_2}{2}-\tfrac{1}{2}+\tfrac{\sqrt{\delta_2}}{2}\right)v_0(x)\right)\phi_R(x)\,dx+L_1+L_2+L_3.
\end{align*}

Using \eqref{assumptions v0,v1 blowup}, we have that there exists $R_1$ such that for any $R\geqslant R_1$
\begin{align*}
\int_{B_R(0)}\left(v_1(x)+\left(\tfrac{\mu_2}{2}-\tfrac{1}{2}+\tfrac{\sqrt{\delta_2}}{2}\right)v_0(x)\right)\phi_R(x)\,dx >0.
\end{align*}

Analogously to the estimates for $K_1,K_2,K_3$, if we introduce
\begin{align*}
\widehat{I}_{\tau,R} & \doteq \int_\frac{\tau}{2}^\tau \int_{B_R} g_1(t) \, \psi_{\tau,R}(t,x)\,  |v(t,x)|^p dx \, dt, \\
\widetilde{I}_{\tau,R} &\doteq \int_0^\tau\int_{B_R(0)\setminus B_{R/2}(0)} g_1(t) \, \psi_{\tau,R}(t,x)\,  |v(t,x)|^p dx \, dt,
\end{align*} then it follows
\begin{align*}
|L_1|+|L_2| &\lesssim \tau^{-2+(\alpha_2-\alpha_1)+\frac{\alpha_1+1}{p'}}R^{\frac{n}{p'}} \widehat{I}_{\tau,R}^{\frac{1}{p}}, \\
|L_3| & \lesssim \tau^{(\alpha_2-\alpha_1)+\frac{\alpha_1+1}{p'}}R^{-2+\frac{n}{p'}} \widetilde{I}_{\tau,R}^{\frac{1}{p}}, 
\end{align*} for $\tau>1$, $R\geqslant R_1$ and provided that
\begin{align*}
(\alpha_2-\alpha_1)p'+\alpha_1>-1 \quad \Leftrightarrow \quad p> \tfrac{1+\alpha_1}{1+\alpha_2}.
\end{align*}

Hence, for $\tau=R> \max\{1,R_1\}$ we obtain
\begin{align} \label{estimate J_R}
J_{R} & \lesssim R^{-2+(\alpha_2-\alpha_1)+\frac{n+\alpha_1+1}{p'}} \Big(\widehat{I}_{\tau,R}^{\frac{1}{p}}+\widetilde{I}_{\tau,R}^{\frac{1}{p}}\Big).
\end{align}

The next step is to combine the estimate for $I_R$ with that one of $J_R$. Of course, $\widehat{J}_R,\widetilde{J}_R\leqslant J_R$, thus, plugging \eqref{estimate J_R} in \eqref{estimate I_R} and conversely, for $\tau=R>\max\{1,R_0,R_1\}$ we have
\begin{align*}
I_R & \lesssim R^{-2+(\alpha_1-\alpha_2)+\frac{n+\alpha_2+1}{q'}} J_R^{\frac{1}{q}} \\
& \lesssim R^{-2+(\alpha_1-\alpha_2)+\frac{n+\alpha_2+1}{q'}+\frac{1}{q}\big(-2+(\alpha_2-\alpha_1)+\frac{n+\alpha_1+1}{p'}\big)} \Big(\widehat{I}_R^{\frac{1}{pq}} +\widetilde{I}_R^{\frac{1}{pq}} \Big), \\
J_{R} & \lesssim R^{-2+(\alpha_2-\alpha_1)+\frac{n+\alpha_1+1}{p'}} I_R^{\frac{1}{p}} \\
& \lesssim R^{-2+(\alpha_2-\alpha_1)+\frac{n+\alpha_1+1}{p'}+\frac{1}{p}\big(-2+(\alpha_1-\alpha_2)+\frac{n+\alpha_2+1}{q'}\big)} \Big(\widehat{J}_R^{\frac{1}{pq}} +\widetilde{J}_R^{\frac{1}{pq}} \Big).
\end{align*}

Let us rewrite the exponents for $R$ in the previous inequalities in a better way. For the first inequality we get
\begin{align*}
-2+(\alpha_1-\alpha_2)+\tfrac{n+\alpha_2+1}{q'}+\tfrac{1}{q}\Big(-2+(\alpha_2-\alpha_1)+\tfrac{n+\alpha_1+1}{p'}\Big)= -2-\tfrac{2}{q}+\big(1-\tfrac{1}{pq}\big) (n+\alpha_1+1)
\end{align*} and for the  second one
\begin{align*}
& -2+(\alpha_2-\alpha_1)+\tfrac{n+\alpha_1+1}{p'}+\tfrac{1}{p}\Big(-2+(\alpha_1-\alpha_2)+\tfrac{n+\alpha_2+1}{q'}\Big)  = -2-\tfrac{2}{p}+\big(1-\tfrac{1}{pq}\big) (n+\alpha_2+1).
\end{align*}

Summarizing, for $\tau=R>\max\{1,R_0,R_1\}$ we have shown
\begin{align}
I_R &  \lesssim R^{-2-\frac{2}{q}+\big(1-\frac{1}{pq}\big) (n+\alpha_1+1)} \Big(\widehat{I}_R^{\frac{1}{pq}} +\widetilde{I}_R^{\frac{1}{pq}} \Big), \label{final estimate IR}\\
J_{R} & \lesssim R^{-2-\frac{2}{p}+\big(1-\frac{1}{pq}\big) (n+\alpha_2+1)} \Big(\widehat{J}_R^{\frac{1}{pq}} +\widetilde{J}_R^{\frac{1}{pq}} \Big). \label{final estimate JR}
\end{align}

Because of the obvious relations $\widehat{I}_R,\widetilde{I}_R\leqslant I_R$, from \eqref{final estimate IR} it follows
\begin{align*}
I_R^{1-\frac{1}{pq}}\lesssim R^{-2(1+\frac{1}{q})+\big(1-\frac{1}{pq}\big) (n+\alpha_1+1)} 
\end{align*} which implies in turn
\begin{align} \label{final estimate IR 2}
I_R \lesssim R^{-2\frac{pq}{pq-1}\big(1+\frac{1}{q}\big)+n+\alpha_1+1} = R^{-2\frac{p(q+1)}{pq-1}+n+\alpha_1+1}.
\end{align}

If the exponent of $R$ on the left hand side is negative, that is, 
\begin{align*}
\tfrac{n+\alpha_1+1}{2} < \tfrac{p(q+1)}{pq-1} \quad\Leftrightarrow \quad \tfrac{n+\alpha_1-1}{2} < \tfrac{p+1}{pq-1}
\end{align*} then, letting $R\to \infty $ in \eqref{final estimate IR 2} we get $\lim_{R\to\infty} I_R=0$. Using the monotone convergence theorem, we find
\begin{align*}
\iint_{[0,\infty)\times \mathbb{R}^n} g_1(t)\,|v(t,x)|^p \, d(t,x)=0
\end{align*} that implies $v=0$ a.e., due to the fact that $g_1$ is always positive. However, this fact contradicts \eqref{assumptions v0,v1 blowup}.

Let us show now that even in the case in which the power of $R$ in \eqref{final estimate IR} is equal to $0$, that is, when $\tfrac{n+\alpha_1-1}{2} = \tfrac{p+1}{pq-1}$, we find the same contradiction. In this last case, \eqref{final estimate IR} implies $I_R\leqslant C$. Hence, by monotone convergence theorem we get
\begin{align*}
\lim_{R\to\infty}I_R= \iint_{[0,\infty)\times \mathbb{R}^n} g_1(t)\,|v(t,x)|^p \, d(t,x) \leqslant C.
\end{align*} So, $g_1|v|^p\in L^1([0,\infty)\times \mathbb{R}^n)$. Consequently, we may employ the dominated convergence theorem, obtaining
\begin{align*}
\lim_{R\to\infty}\widehat{I}_R=0 \quad \mbox{and} \quad \lim_{R\to\infty}\widetilde{I}_R=0.
\end{align*} Using these relations in \eqref{final estimate IR}, we find as in the previous case $\lim_{R\to\infty} I_R=0$. Repeating the previous argument, we arrive at the same contradiction.

In an analogous way, one can show that \eqref{final estimate JR} leads to the condition $u=0$ a.e. in the case in which $\frac{n+\alpha_2-1}{2} \leqslant \frac{q+1}{pq-1}$, but this fact is not possible because of \eqref{assumptions u0,u1 blowup}. Summarizing, we proved that for 
\begin{align*}
\tfrac{n+\alpha_1-1}{2} \leqslant \tfrac{p+1}{pq-1} \quad \mbox{or} \quad \tfrac{n+\alpha_2-1}{2} \leqslant \tfrac{q+1}{pq-1},
\end{align*}  provided that $p,q$ fulfill \eqref{lower bound p or q blowup}, the weak solution $(u,v)$ cannot be globally in time defined. Since the first previous relations on $(p,q)$ are equivalent to \eqref{condition p,q blowup}, the proof is completed.

\section{Proofs of global existence results} \label{Section Proofs SDGE results}

Let us introduce some common notations for the proofs of the global (in time) existence results. We denote by $E^{(\mu,\nu)}_0(t,s,x)$ and $E^{(\mu,\nu)}_1(t,s,x)$ the fundamental solutions to \eqref{scale inv eq lin}, that is, the distributional solutions to \eqref{scale inv eq lin} with initial data $(u_0,u_1)=(\delta_0,0)$ and $(u_0,u_1)=(0,\delta_0)$, respectively. Hence, the solution to \eqref{scale inv eq lin} is given by
\begin{align*}
u(t,x)=E^{(\mu,\nu)}_0(t,s,x) \ast_{(x)} u_0(x)+E^{(\mu,\nu)}_1(t,s,x) \ast_{(x)} u_1(x).
\end{align*}

Having in mind Duhamel's principle, let us introduce the operator
\begin{align*}
N: (u,v) \to N(u,v)\doteq \Big( u^{\lin}+G_1(v),v^{\lin}+G_2(u)\Big),
\end{align*} where $(u^{\lin},v^{\lin})$ is the solution to the corresponding linear homogeneous system with data $(u_0,u_1;v_0,v_1)$, that is,
\begin{align*}
u^{\lin}(t,x) & \doteq E^{(\mu_1,\nu_1)}_0(t,0,x) \ast_{(x)} u_0(x)+E^{(\mu_1,\nu_1)}_1(t,0,x) \ast_{(x)} u_1(x), \\
v^{\lin}(t,x) & \doteq E^{(\mu_2,\nu_2)}_0(t,0,x) \ast_{(x)} v_0(x)+E^{(\mu_2,\nu_2)}_1(t,0,x) \ast_{(x)} v_1(x), 
\end{align*} and $G_1(v),G_2(u)$ are the following integral operators:
\begin{align*}
G_1(v)(t,x) & \doteq \int_0^t E^{(\mu_1,\nu_1)}_1(t,s,x) \ast_{(x)} |v(s,x)|^p \, ds, \\
G_2(u)(t,x) & \doteq \int_0^t E^{(\mu_2,\nu_2)}_1(t,s,x) \ast_{(x)} |u(s,x)|^q \, ds.
\end{align*}

Moreover, we introduce a family of function spaces $\{X(T)\}_{T>0}$, with 
\begin{align} \label{def X(T)}
X(T) \doteq \big(\mathcal{C}([0,T],H^1(\mathbb{R}^n))\cap\mathcal{C}^1([0,T],L^2(\mathbb{R}^n))\big)^2
\end{align} equipped with the norm 
\begin{align}
\|(u,v)\|_{X(T)}\doteq \sup_{t\in[0,T]} \Big((1+t)^{-\gamma_1}M_1(t,u)+(1+t)^{-\gamma_2}M_2(t,v)\Big),\label{def norm X(T)}
\end{align} where
\begin{align*}
M_1(t,u) &\doteq (1+t)^{\frac{n}{2}+\alpha_1}\Big(\|(\nabla u,u_t)(t,\cdot)\|_{L^2(\mathbb{R}^n)}+(1+t)^{-1}\|u(t,\cdot)\|_{L^2(\mathbb{R}^n)}\Big), \\
M_2(t,v) &\doteq (1+t)^{\frac{n}{2}+\alpha_2}\Big(\|(\nabla v,v_t)(t,\cdot)\|_{L^2(\mathbb{R}^n)}+(1+t)^{-1}\|v(t,\cdot)\|_{L^2(\mathbb{R}^n)}\Big),
\end{align*}
and  $\gamma_1,\gamma_2\geqslant 0$ represent possible losses of decay for $(u,v)$ in comparison with the corresponding decay estimates for $(u^{\lin},v^{\lin})$.

In order to prove the global (in time) existence of solutions to \eqref{weakly coupled system} we want to prove that the operator $N$ is a contraction on $X(T)$ with an independent of $T$ Lipschitz constant. Then, the solution $(u,v)$ to \eqref{weakly coupled system} will be the solution of the nonlinear integral system of equation $(u,v)=N(u,v)$, i.e., the unique fixed point of $N$. More specifically, we will prove the inequalities
\begin{align}
\| N(u,v)\|_{X(T)} &\lesssim \|(u_0,u_1)\|_{\mathcal{A}}+\|(v_0,v_1)\|_{\mathcal{A}} + \| (u,v)\|_{X(T)}^p+ \| (u,v)\|_{X(T)}^q \label{1st ineq contr}\\
\| N(u,v)-N(\bar{u},\bar{v})\|_{X(T)} &\lesssim 
\| (u,v)-(\bar{u},\bar{v})\|_{X(T)} \bigg(\sum_{r=p,q}\| (u,v)\|_{X(T)}^{r-1}+\| (\bar{u},\bar{v})\|_{X(T)}^{r-1}\bigg)  \label{2nd ineq contr}
\end{align}
uniformly with respect to $T$, which imply the desired property for the operator $N$, provided that $\|(u_0,u_1)\|_{\mathcal{A}}+\|(v_0,v_1)\|_{\mathcal{A}}\doteq \varepsilon$ is sufficiently small. 

Let us underline explicitly, that \eqref{1st ineq contr} and \eqref{2nd ineq contr} imply for the fixed point $(u,v)$ of $N$ the estimates
\begin{align*}
M_1(t,u) & \lesssim \varepsilon\, (1+t)^{\gamma_1}, \\
M_2(t,v) & \lesssim \varepsilon \,(1+t)^{\gamma_2}, 
\end{align*} which are exactly the estimates for $(u,v)$ in Theorems \ref{Thm global existence supercritical p,q}, \ref{Thm global existence subcritical p and supercritical q} and \ref{Thm global existence supercritical p and subcritical q} provided that $\gamma_1$ and $\gamma_2$ are suitably choosen (for example, at least one among them has to be 0).

\subsection{Proof of Theorem \ref{Thm global existence supercritical p,q}} \label{Section supercritical p,q}

Let us consider the space $X(T)$ defined by \eqref{def X(T)} and equipped with the norm given by \eqref{def norm X(T)} with $\gamma_1=\gamma_2=0$. Due to the fact that we are in the subcase \eqref{supercritical p and q}, no loss of decay is required in comparison to the homogeneous linear problem neither for $u$ nor for $v$. From Proposition \ref{Prop Lin Estim} it follows immediately
\begin{align*}
\| (u^{\lin},v^{\lin})\|_{X(T)}\lesssim \|(u_0,u_1)\|_{\mathcal{A}}+\|(v_0,v_1)\|_{\mathcal{A}}.
\end{align*} Consequently, in order to show \eqref{1st ineq contr} it remains to prove that
\begin{align}\label{3rd ineq contr}
\| (G_1(v),G_2(u))\|_{X(T)}\lesssim  \| (u,v)\|_{X(T)}^p+ \| (u,v)\|_{X(T)}^q.
\end{align}

Let us begin by estimating $M_1(t,G_1(v))$. For $j+\ell=0,1$, by Proposition \ref{Prop Lin Estim v_0=0 for t>s} we have
\begin{align*}
\| \nabla^j \partial_t^\ell G_1(v)(t,\cdot)& \|_{L^2(\mathbb{R}^n)} \leqslant \int_0^t \| E^{(\mu_1,\nu_1)}(t,s,x) \ast_{(x)} |v(s,\cdot)|^p \|_{L^2(\mathbb{R}^n)} \,ds \\
&\lesssim (1+t)^{-(j+\ell)-\frac{n}{2}-\alpha_1+1} \int_0^t (1+s)^{\alpha_1} \Big(\| v(s,\cdot) \|_{L^p(\mathbb{R}^n)}^p+(1+s)^{\frac{n}{2}}\| v(s,\cdot) \|_{L^{2p}(\mathbb{R}^n)}^p\Big) \,ds.
\end{align*}
Using Gagliardo-Nirenberg inequality, we can estimate the $L^p$ norm and the $L^{2p}$ norm of $v(s,\cdot)$ as follows:
\begin{align*}
\| v(s,\cdot)\|_{L^{hp}(\mathbb{R}^n)} & \lesssim \| v(s,\cdot)\|_{L^2(\mathbb{R}^n)}^{1-\theta(hp)} \| \nabla v(s,\cdot)\|_{L^2(\mathbb{R}^n)}^{\theta(hp)} \\
& \lesssim (1+s)^{-\frac{n}{2}-\alpha_1+1-\theta(hp)} \| (u,v)\|_{X(s)},
\end{align*} for $h=1,2$, where $\theta(r)\doteq n\big(\frac{1}{2}-\frac{1}{r}\big)\in [0,1]$. In particular, the belonging of $\theta(p)$ and $\theta(2p)$ to the interval $[0,1]$ implies the further conditions $2\leqslant p$ and $p\leqslant \frac{n}{n-2}$ for $n\geqslant 3$. Also, 
\begin{align}
(1+t)^{(j+\ell)+\frac{n}{2}+\alpha_1-1} \| \nabla^j \partial_t^\ell G_1(v)(t,\cdot)\|_{L^2(\mathbb{R}^n)}  &\lesssim  \int_0^t (1+s)^{-(n+\alpha_2-1)p+n+\alpha_1}  \| (u,v)\|_{X(s)}^p  \,ds  \notag \\&\lesssim    \| (u,v)\|_{X(t)}^p, \label{integral |v|^p}
\end{align} where we used the condition $p>\widetilde{p}\,(n,\alpha_1,\alpha_2)$ in order to guarantee the uniform boundedness of the integral in the last inequality.

Similarly, we can estimate $M_2(t,G_2(u))$ in the following way
\begin{align}
\| \nabla^j \partial_t^\ell G_2(u)(t,\cdot)&\|_{L^2(\mathbb{R}^n)}\leqslant \int_0^t \| E^{(\mu_2,\nu_2)}(t,s,x) \ast_{(x)} |u(s,\cdot)|^q \|_{L^2(\mathbb{R}^n)} \,ds  \notag \\
&\lesssim (1+t)^{-(j+\ell)-\frac{n}{2}-\alpha_2+1} \int_0^t (1+s)^{\alpha_2} \Big(\| u(s,\cdot) \|_{L^q(\mathbb{R}^n)}^q+(1+s)^{\frac{n}{2}}\| u(s,\cdot) \|_{L^{2q}(\mathbb{R}^n)}^q\Big) \,ds  \notag \\
&\lesssim (1+t)^{-(j+\ell)-\frac{n}{2}-\alpha_2+1} \int_0^t (1+s)^{-(n+\alpha_1-1)q+n+\alpha_2}  \| (u,v)\|_{X(s)}^q  \,ds  \notag  \\&\lesssim    (1+t)^{-(j+\ell)-\frac{n}{2}-\alpha_2+1} \| (u,v)\|_{X(t)}^q \label{integral |u|^q}
\end{align} for $j+\ell=0,1$, where in the last step we employed the assumption $q> \widetilde{q}\,(n,\alpha_1,\alpha_2)$. Combining \eqref{integral |v|^p} and \eqref{integral |u|^q} we get immediately \eqref{3rd ineq contr}.

Let us sketch briefly the proof of the Lipschitz condition \eqref{2nd ineq contr}. 
As $$N(u,v)-N(\bar{u},\bar{v})=\big(G_1(v)-G_1(\bar{v}),G_2(u)-G_2(\bar{u})\big),$$ it is sufficient to control the quantities $M_1(t,G_1(v)-G_1(\bar{v}))$ and $M_2(t,G_2(u)-G_2(\bar{u}))$.
Using again Proposition \ref{Prop Lin Estim v_0=0 for t>s}, we find
\begin{align*}
& (1+t)^{(j+\ell)+\frac{n}{2}+\alpha_1-1} \| \nabla^j \partial_t^\ell (G_1(v)-G_1(\bar{v}))(t,\cdot)\|_{L^2(\mathbb{R}^n)} \\
& \qquad \lesssim  \int_0^t (1+s)^{\alpha_1} \Big(\| |v(s,\cdot)|^p- |\bar{v}(s,\cdot)|^p\|_{L^1(\mathbb{R}^n)}+(1+s)^{\frac{n}{2}}\| |v(s,\cdot)|^p- |\bar{v}(s,\cdot)|^p \|_{L^{2}(\mathbb{R}^n)}\Big) \,ds.
\end{align*}

By the pointwise estimate $| |v|^p-|\bar{v}|^p| \leqslant p (|v|^{p-1}+|\bar{v}|^{p-1})|v-\bar{v}|$, H\"older's inequality, Gagliardo-Nirenberg inequality and the definition of norm for the family of spaces $\{X(t)\}_{t>0}$, for $h=1,2$ we arrive at 
\begin{align*}
\| |v(s,\cdot)|^p- |\bar{v}(s,\cdot)|^p\|_{L^h(\mathbb{R}^n)} & \lesssim \| v(s,\cdot)- \bar{v}(s,\cdot)\|_{L^{ph}(\mathbb{R}^n)}  \Big(\| v(s,\cdot)\|_{L^{hp}(\mathbb{R}^n)}^{p-1} +\| \bar{v}(s,\cdot)\|_{L^{hp}(\mathbb{R}^n)}^{p-1} \Big) \\
& \lesssim (1+s)^{-(n+\alpha_2-1)+\frac{n}{h}}\| (u,v)- (\bar{u},\bar{v})\|_{X(s)}  \Big(\| (u,v)\|_{X(s)}^{p-1} +\| (\bar{u},\bar{v})\|_{X(s)}^{p-1} \Big).
\end{align*} So, combining the last two estimates, we get 
\begin{align*}
 (1+t)^{(j+\ell)+\frac{n}{2}+\alpha_1-1} \| \nabla^j \partial_t^\ell (G_1(v)- & G_1(\bar{v}))(t,\cdot)\|_{L^2(\mathbb{R}^n)} \\ & \lesssim \| (u,v)- (\bar{u},\bar{v})\|_{X(t)}  \Big(\| (u,v)\|_{X(t)}^{p-1} +\| (\bar{u},\bar{v})\|_{X(t)}^{p-1} \Big),
\end{align*} provided that $p>\widetilde{p}\,(n,\alpha_1,\alpha_2)$. In an analogous way, we can prove \begin{align*}
 (1+t)^{(j+\ell)+\frac{n}{2}+\alpha_2-1} \| \nabla^j \partial_t^\ell (G_2(u)- & G_2(\bar{u}))(t,\cdot)\|_{L^2(\mathbb{R}^n)} \\ &\lesssim \| (u,v)- (\bar{u},\bar{v})\|_{X(t)}  \Big(\| (u,v)\|_{X(t)}^{q-1} +\| (\bar{u},\bar{v})\|_{X(t)}^{q-1} \Big),
\end{align*} due to $q>\widetilde{q}\,(n,\alpha_1,\alpha_2)$. Thus, we proved \eqref{2nd ineq contr}. This concludes the proof.

\subsection{Proof of Theorem \ref{Thm global existence subcritical p and supercritical q}} \label{Section subcritical p, supercritical q}

Let us consider the space $X(T)$ defined by \eqref{def X(T)} and equipped with the norm given by \eqref{def norm X(T)} with $$\gamma_1\equiv\gamma = \begin{cases}(n+\alpha_2-1)(\widetilde{p}\,(n,\alpha_1,\alpha_2)-p) & \mbox{if} \ \ p < \widetilde{p}\,(n,\alpha_1,\alpha_2), \\ \epsilon & \mbox{if} \ \ p = \widetilde{p}\,(n,\alpha_1,\alpha_2), \end{cases}$$ defined as in the statement and $\gamma_2=0$. Since we assume $p\leq \widetilde{p}\,(n,\alpha_1,\alpha_2)$, in order to control the integral in \eqref{integral |v|^p} we allow this loss of decay for the estimates for $u$ in comparison with those for $u^{\lin}$.  By Proposition \ref{Prop Lin Estim} we get
\begin{align*}
\| (u^{\lin},v^{\lin})\|_{X(T)}\lesssim \sup_{t\in[0,T]}(1+t)^{-\gamma} \|(u_0,u_1)\|_{\mathcal{A}}+\|(v_0,v_1)\|_{\mathcal{A}} \lesssim \|(u_0,u_1)\|_{\mathcal{A}}+\|(v_0,v_1)\|_{\mathcal{A}},
\end{align*} due to the fact that $\gamma>0$. Let us consider now the nonlinear part $(G_1(v),G_2(u))$. As in the previous section, we have
\begin{align*}
(1+t)^{(j+\ell)+\frac{n}{2}+\alpha_1-1} \| \nabla^j \partial_t^\ell G_1(v)(t,\cdot)\|_{L^2(\mathbb{R}^n)}  &\lesssim  \int_0^t (1+s)^{-(n+\alpha_2-1)p+n+\alpha_1}  \| (u,v)\|_{X(s)}^p  \,ds 
\end{align*} but, differently from the previous section we can no longer estimate the last integral uniformly by a constant, as we are in the case $p\leqslant \widetilde{p}\,(n,\alpha_1,\alpha_2)$. So, we have
\begin{align*}
\int_0^t (1+s)^{-(n+\alpha_2-1)p+n+\alpha_1} \lesssim (1+t)^{\gamma},
\end{align*} which implies the desired estimate
\begin{align*}
(1+t)^{(j+\ell)+\frac{n}{2}+\alpha_1-1-\gamma} \| \nabla^j \partial_t^\ell G_1(v)(t,\cdot)\|_{L^2(\mathbb{R}^n)}  &\lesssim  \| (u,v)\|_{X(t)}^p.
\end{align*} The next step is to determine under which condition for $(p,q)$ the inequality $M_2(t,G_2(u)) \lesssim \| (u,v)\|_{X(t)}^q$ holds. Similarly to the previous section, keeping in mind that now we have a different decay rate for $u$ coming for the norm of $(u,v)\in X(s)$, we get 
\begin{align*}
\| \nabla^j \partial_t^\ell G_2(u)(t,\cdot)& \|_{L^2(\mathbb{R}^n)} 
\\ &\lesssim (1+t)^{-(j+\ell)-\frac{n}{2}-\alpha_2+1} \int_0^t (1+s)^{\alpha_2} \Big(\| u(s,\cdot) \|_{L^q(\mathbb{R}^n)}^q+(1+s)^{\frac{n}{2}}\| u(s,\cdot) \|_{L^{2q}(\mathbb{R}^n)}^q\Big) \,ds  \\
&\lesssim (1+t)^{-(j+\ell)-\frac{n}{2}-\alpha_2+1} \int_0^t (1+s)^{(1+\gamma-n-\alpha_1)q+n+\alpha_2}    \,ds \,   \| (u,v)\|_{X(t)}^q \\&\lesssim    (1+t)^{-(j+\ell)-\frac{n}{2}-\alpha_2+1} \| (u,v)\|_{X(t)}^q, 
\end{align*} provided that the exponent of the integrand is smaller than $-1$. If $p<\widetilde{p}(n,\alpha_1,\alpha_2)$, this condition is equivalent to require
\begin{align*}
q(1+(n+\alpha_1+1)-(n+\alpha_2-1)p-n-\alpha_1)+ & n+\alpha_2+1<0 \\ & \quad \Leftrightarrow \quad 2(q+1)+(n+\alpha_2-1)(1-pq)<0 \\
& \quad \Leftrightarrow \quad \tfrac{q+1}{pq-1}<\tfrac{n+\alpha_2-1}{2}
\end{align*} that is, for $(p,q)$ fulfilling \eqref{F(q,p,alpha2)<0}. On the other hand, for $p=\widetilde{p}(n,\alpha_1,\alpha_2)$, we can choose $\epsilon>0$ so small that
\begin{align*}
q(2+\epsilon-(n+\alpha_2-1)p)+n+\alpha_2+1<0 & \quad \Leftrightarrow \quad 
\tfrac{\big(1+\tfrac{\epsilon}{2}\big)q+1}{pq-1}<\tfrac{n+\alpha_2-1}{2}
\end{align*} is satisfied. In both cases \eqref{F(q,p,alpha2)<0} implies the desired inequality. Hence, we proved \eqref{1st ineq contr}. The proof of \eqref{2nd ineq contr} is completely similar to that one in the proof of Theorem \ref{Thm global existence supercritical p,q}. So, the proof is over.

\begin{rem} It is clear that, due to the symmetry of \eqref{weakly coupled system}, the proof of Theorem \ref{Thm global existence supercritical p and subcritical q} is completely analogous to that one of Theorem \ref{Thm global existence subcritical p and supercritical q} by choosing $\gamma_1=0$ and $\gamma_2\equiv \bar{\gamma}$ as in the statement of Theorem \ref{Thm global existence supercritical p and subcritical q}.
\end{rem}

\section{Final remarks}

Combining the results from Theorems \ref{Thm blow-up MTF}, \ref{Thm global existence supercritical p,q}, \ref{Thm global existence subcritical p and supercritical q} and \ref{Thm global existence supercritical p and subcritical q} we have that $$\max \big\{\tfrac{p+1}{pq-1}-\tfrac{\alpha_1}{2}, \tfrac{q+1}{pq-1}-\tfrac{\alpha_2}{2}\big\} =\tfrac{n-1}{2}$$ is the critical exponent for the weakly coupled system \eqref{weakly coupled system} provided that the coefficients satisfy $\delta_1,\delta_2>(n+1)^2$ in the sense we explained in Section \ref{Section main results}. Actually, one can slightly improve this result up to range $\delta_1,\delta_2\geqslant (n+1)^2$ modulo a (possible) further arbitrarily small loss of decay rate with respect to the case $\delta_1,\delta_2>(n+1)^2$  in Theorems \ref{Thm global existence subcritical p and supercritical q} and \ref{Thm global existence supercritical p and subcritical q}

In the case $0<\delta_1< (n+1)^2$ or $0<\delta_2< (n+1)^2$ we cannot obtain a sharp result as in the above mentioned case by using $L^2-L^2$ estimates with additional $L^1$ regularity and working in  classical energy spaces, due to the fact that the first order derivatives have a weaker decay rate (cf. Theorems 4.6 and 4.7 in \cite{PalRei17} for further details).

In the case in which $\mu_1=\mu_2$ and $\nu_1=\nu_2$, the critical exponent for \eqref{weakly coupled system} is
\begin{align*}
\tfrac{\max\{p,q\}}{pq-1}-\tfrac{n+\alpha-1}{2}=0,
\end{align*} where $\alpha= \alpha_1=\alpha_2$. In particular,
\begin{align*}
& \tfrac{\max\{p,q\}}{pq-1}-\tfrac{n+\alpha-1}{2}\geqslant 0 \qquad \Leftrightarrow \qquad  \max\{p,q\} \big(\min\{p,q\}+1-p_{\Fuj}(n+\alpha-1)\big) \leqslant p_{\Fuj}(n+\alpha-1), \\
& \tfrac{\max\{p,q\}}{pq-1}-\tfrac{n+\alpha-1}{2}< 0 \qquad \Leftrightarrow \qquad  \max\{p,q\} \big(\min\{p,q\}+1-p_{\Fuj}(n+\alpha-1)\big) > p_{\Fuj}(n+\alpha-1),
\end{align*} where $p_{\Fuj}(n+\alpha-1)$ is the critical exponent for the corresponding single equation (see also \cite{NPR16,PalRei17,Pal17}).

Since for the single equation \eqref{scale inv eq} we expect $p_0(n+\mu)$ to be the critical exponent for small and nonnegative values of $\delta$ (cf. \cite{NPR16,PR17vs,Pal18odd,Pal18even,PT18}), it is clear that the result from Theorem \ref{Thm blow-up MTF} cannot be sharp in this case. 

Indeed, in an upcoming paper a blow-up result for \eqref{weakly coupled system} is going to be proved for $\delta_1,\delta_1\geqslant 0$ provided that $p,q>1$ satisfy $$\max \Big\{\frac{p+2+q^{-1}}{pq-1}-\frac{\mu_1}{2}, \frac{q+2+p^{-1}}{pq-1}-\frac{\mu_2}{2}\Big\} >\frac{n-1}{2}.$$ We notice that the corresponding critical relation for the pair $(p,q)$ is a shift of the critical exponent for \eqref{weaklycoupledwave}.

\section*{Acknowledgments} 

The PhD study of the first author is supported by S\"achsiches Landesgraduiertenstipendium.
The second author is member of the Gruppo Nazionale per L'Analisi Matematica, la Probabilit\`{a} e le loro Applicazioni (GNAMPA) of the Instituto Nazionale di Alta Matematica (INdAM). The authors thank their supervisor Michael Reissig for the suggestions in the preparation of the final version.

\bibliographystyle{elsarticle-num}

\end{document}